\documentclass{amsproc}
\usepackage{amssymb, latexsym}
\newcommand{\R}{{\mathbb R}}

\newcommand{\T}{{\mathbb T}}

\newcommand{\N}{{\mathbb N}}

\def\cal{\mathcal}

\newcommand{\p}{\partial}
\newcommand{\Ly}{{\cal{L}}}

\newtheorem{theorem}{Theorem}[section]
\theoremstyle{definition}

\theoremstyle{remark}

\theoremstyle{proposition}
\newtheorem{prop}[theorem]{Proposition}
\theoremstyle{lemma}
\newtheorem{lemma}[theorem]{Lemma}
\theoremstyle{corollary}

\theoremstyle{claim}
\newtheorem{claim}{Claim}

\newenvironment{Proof}{\noindent {\sc Proof}}{$\Box$ \vspace{2 ex}}

\numberwithin{equation}{section}

\begin{document}

\title{ Stability of Solitons for the KdV equation in $H^s$, $0 \leq s
< 1$ \\
\emph{Preliminary Version}}

\author{S.~Raynor and G.~Staffilani}
\address{Massachusetts Institute of Technology}
\thanks{G.S. is supported in part by N.S.F. Grant DMS 0100375
and a grant by the Sloan Foundation.}
\email{sraynor@math.mit.edu and gigliola@math.mit.edu}


\begin{abstract}
We study the long-time stability of soliton solutions to the
Korteweg-deVries equation.  We consider solutions $u$ to the KdV with
initial data in $H^s$, $0 \leq s < 1$, that are initially close in
$H^s$ norm to a soliton.   We prove that the possible orbital
instability of these ground states is at most polynomial in time.
This is an analogue to the $H^s$ orbital instability result of
\cite{CKSTT3}, and obtains the same maximal growth rate in $t$.  Our
argument is based on the {\lq\lq}$I$-method{\rq\rq} used in \cite{CKSTT3} and
other papers of Colliander, Keel, Staffilani, Takaoka and Tao, which
pushes these $H^s$ functions to the $H^1$ norm.  
\end{abstract}

\maketitle

\section{Introduction}

We will consider the long-time stability of soliton solutions to the
Korteweg-deVries Equation.  The KdV equation, which was developed as a
model for one-dimensional waves in shallow water, is as follows:
\begin{equation}\label{kdv}
u_t + u_{xxx} + (u^2)_x = 0.
\end{equation}
We will consider the initial value problem for the KdV with initial
data $u_0 \in H^s, 0 \leq s < 1$.  Local well-posedness (that is,
short-time existence, uniqueness and uniform continuity with regard to
initial data) for the Cauchy problem is known. (See \cite{B} and \cite{KPV} for the most recent results.)
Moreover, the KdV equation has an infinite sequence of conservation
laws which hold for any solution which is sufficiently smooth.  The
first few are:
\begin{align*}
 G(u) & = \int_\R u(x,t) dx = \int_\R u(x,0) dx,  \\
\| u (t) \|^2_{L^2} &= \int_\R |u(x,t)|^2 dx = \int_\R |u(x,0)|^2 dx,  \\
H(u) &= \int_\R ( |\p _x u(x,t)|^2 - \frac23 u(x,t)^3) dx = \int_\R ( |\p
_x u(x,0)|^2 - \frac23 u(x,0)^3) dx.  
\end{align*}
Using the local well-posedness arguments, these
conservation laws, and iteration arguments, global well-posedness can be 
deduced for $s \geq 0$.\footnote{Global well-posedness also holds for
$s > -\frac34$ \cite{CKSTT4}.}   

It is known that the KdV equation
admits traveling wave solutions called solitons which satisfy $ Q(x,
t) = \psi(x - Ct)$, and $\psi$ therefore is a  solution to the following
ODE:
\begin{equation}\label{ODE}
\psi_{xx} - C \psi + \psi^2 = 0.
\end{equation}
There exists a unique even, positive solution $\psi_0$ to this equation.  This
soliton is smooth and rapidly decreasing as $|x| \rightarrow 
\infty$.  In fact, 
\begin{equation}\label{psi}
\psi_0(x) = \frac32 C \mathrm{sech} ^2 ( \frac12 C^\frac12 x ).
\end{equation}
For simplicity we will consider only the case $C = 1$, as the others
can be recovered by scaling.  We will define $\Sigma = \{ \psi_0 (x - x_0) |
x_0 \in \R \}$ to be the one-parameter space of all solitons moving
with speed $1$.  Note that the KdV flow preserves $\Sigma$ and that
each element of $\Sigma$ is a solution to (\ref{ODE}).  

It was proven by Benjamin \cite{Be} in 1972 that soliton solutions are
stable in the following sense: if $u$ is a solution to the KdV which
is initially close to a soliton in $H^1$ norm, then for all time $u$
is close to a soliton.  Some corrections and extensions of his result
were offered by Bona \cite{Bo}.  More recently, Weinstein \cite{W} has
offered a general theory which proves the stability of soliton
solutions to generalized KdV equations as well as a class of
non-linear Schr$\mathrm{\ddot{o}}$dinger equations.  In \cite{CKSTT3},
Colliander, Keel, Staffilani, Takaoka, and Tao exploited Weinstein's
result to prove that the instability of soliton solutions to the NLS
in $H^s, 0 \leq s < 1$ grows at most polynomially in $t$.  They made
use of a multiplier operator which they had developed in their proof
of global well-posedness for dispersive equations with initial data in
$H^s$, $0 \leq s < 1$.\cite{CKSTT4} This multiplier operator allowed
them to work with $H^1$ norms, which they could then control using
Weinstein's result.

In this paper, we will again exploit the multiplier operator which
they developed, as well as the original proof of $H^1$ stability of
solitons for the KdV.  We will prove that in $H^s, 0 \leq s < 1$,
soliton solutions to the KdV are at most polynomially unstable.  Our
main result is:

\begin{theorem}\label{main}    

Let $0 \leq s < 1$, Let $\sigma = \mathrm{dist}_{H^s} (u_0, \Sigma) \ll 1$, and
let $u$ be the solution to the KdV such that $u ( \cdot , 0) = u_0$.
Then $\mathrm{dist}_{H^s} (u(t), \Sigma) \leq t^{1-s+\epsilon}
\sigma$, for all $t$ such that $t \ll \sigma^{-\frac1{1-s+\epsilon}}.$

\end{theorem}

To prove this, we will employ the Lyapunov functional introduced by
Benjamin \cite{Be}:  
$$\Ly (u) = \|u\|^2_{L^2} + H(u) = \int_\R |u_x|^2 + |u|^2
-\frac23|u|^3 .$$  It can be shown using the Gagliardo-Niremberg
inequality that $\Ly \geq 0.$  Note that if $u$ is a solution to the KdV equation with $u \in H^s$
and $s \geq 1$, then $\Ly(u)$ is conserved.  In fact, we have the
equation
$$\p _t \Ly (u) = 2 \int \limits _\R u_t ( -u_{xx} + u - u^2) dx , $$
which vanishes if $u$ is a solution to (\ref{kdv}) by integration by
parts.  This calculation also shows that solitons, which are
solutions to (\ref{ODE}), are critical points of the functional $\Ly$.
In \cite{Be} (see also \cite{W}), Benjamin proved that they are
minimizers and moreover 
that, for all $u \in H^1$ such that $\mathrm{dist} _{H^1} (u, \Sigma) \ll 1$, 
\begin{equation}
\Ly(u) - \Ly(Q) \sim \mathrm{dist} _{H^1} (u, \Sigma) ^2 . \label{H1}
\end{equation}  
This then implies the stability of the solitons because $\Ly (u)$ is
conserved in $t$.

We will extend this result to $H^s$, $0 \leq s < 1$, finding the
possible growth in time of the distance between $u$ and the solitons
to be at most polynomial.  To do so, we will exploit the fact that the quantity
$\Ly(Iu)$ is almost conserved in time, where $I$ is a smoothing
operator that maps $H^s$ to $H^1$.  This techniques was used by
Colliander, Keel, Staffilani, Takaoka, and Tao in \cite{CKSTT3} to
prove polynomial stability bounds for solitons solutions to the
Schr$\mathrm{\ddot{o}}$dinger Equation.  We will follow the technique
developed in that paper in general outline, making the necessary
estimates for the KdV equation.  We will also follow the structure of
that paper, giving progressively more sophisticated arguments that get
closer to Theorem \ref{main} with each iteration.  

Several interesting open questions remain.  It is not known whether
the power of $t$ which we obtain in the theorem is sharp.  Moreover,
we have not completed the estimates for the modified KdV equation
$$u_t + u_{xxx} + (u^3)_x = 0$$
and it is not known whether such stability results hold in that case.
Finally, a recent paper of Merle and Vega \cite{MV} has concluded that
in fact KdV solitons are stable in $L^2$.\footnote{The reader may
think that some sort of ``interpolation'' between the $H^1$ and $L^2$
stability should give an even better result than the authors obtain,
but unfortuntely it is no obvious how to ``interpolate.''}  The
authors are currently studying whether this 
result and the $I$-method exploited in this paper can be extended to
prove polynomial stability bounds below $L^2$.  

The structure of this paper is as follows.  In section 2, we define
our notation and quote some important estimates that will be used in
the following sections.  In section 3, we make our first attempt at
proving the main theorem, obtaining a weaker form of the estimate.  In
section 4, we refine the techniques of section 3 but still miss the
main theorem by an $\epsilon$ power in $\mathrm{dist}_{H^s} (u_0,
\Sigma)$.  Finally, in section 5 we complete the proof of the
theorem.  

\section{Notation and Set-Up}

We will use the notation $A \lesssim B$ to mean that $A \leq c B$
where $c$ is a constant depending on $s$ that may vary from line to
line, and similarly for the notation $A \sim B$.  We will use $\langle
\xi \rangle$ to denote $1 + |\xi|$.   

We define the spatial Fourier transform by 
$$\hat{f}(\xi) = \int \limits_\R e^{-i x \xi} f(x) dx$$
and the spacetime Fourier transform by
$$\tilde{u}(\xi, \tau) = \int \limits_\R e^{-i ( x \xi + t \tau)} f(x,t) dxdt.$$
We define the $X^{s,b}$ space, as in \cite{B}, by the norm 
$$ \| u \|_{X^{s,b}} = \| \langle \xi \rangle^s \langle \tau - \xi ^3
\rangle ^b \tilde{u} (\xi, \tau) \|_{L^2_{\xi, \tau}}.$$
We will also use the notation 
$$X^{s,b} _I = \{ u \arrowvert _ {\R \times I} \colon u \in X^{s, b} \} $$
with the norm
$$\| u \|_{X^{s,b}_I} = \inf \{ \|v \|_{X^{s, b}} \colon v \arrowvert _{\R
\times I} = u \}.$$
We will use the notation
$$\Lambda _n(m(\xi_1, \ldots, \xi_n); f_1, \ldots, f_n) =
\hskip -.1in 
\int \limits _{\xi_1 + \ldots + \xi_n = 0 } \hskip -.1in [m(\xi_1,
\ldots, \xi_n)] \hat{f_1}(\xi_1) \cdot \ldots \cdot
\hat{f_n}(\xi_n) d\xi_1 \ldots d\xi_n,$$
where $[m(\xi_1, \ldots, \xi_n)]$ is the symmetrization of the
multiplier $m$ in the $\xi_i$ variables.  Note that we will not always
work with the symmetrized multiplier if it does not matter, but that
occasionally symmetrization will be necessary to obtain the
appropriate estimates.  

For $N \gg 1$ and fixed, we define the operator $I = I _N$ to be a
smooth even multiplier operator such that 
\begin{equation}\label{I}
\widehat{(I_N u)}= \begin{cases} \hat{u} & |\xi| < N \\
                         \frac{\xi^{s-1}}{N^{s-1}} \hat{u} & |\xi| > 10N 
                \end{cases}.
\end{equation}
We generally omit the subscript $N$ unless it is necessary for clarity. 
We also use the notation $N_i$ for a dyadic block in 
the frequency space of the function $u_i$, that is, in the domain of the variable $\xi_i$.  Note that $N_i$ is not necessarily positive.  We will write $u_{i,N_i}$ for
the function obtained from $u_i$ by restricting it to its components
with frequency in $N_i$.  That is, if $\phi_{N_i}$ is a smooth cutoff
function which is the identity in $[N_i, 2N_i]$ and which has support in $[N_i
-1, 2N_i + 1]$, then $\widehat{u_{i, N_i}} = \phi_{N_i} \widehat{u_i}$.  

Finally, we will denote by $W(t)$ the solution operator for the linear
KdV equation, $u_t + u_{xxx} = 0$.  

We will also use the following estimates:

\begin{enumerate}
\item By the Plancherel Theorem and Cauchy-Schwartz, we have\footnote{For a more precise proof of this estimate on a
finite interval in $t$, see \cite{CKSTT3a}} and the proof of Lemma 8.1
in \cite{CKSTT3b}:
\begin{equation}
\int \int |{u_1} {u_2} | \ dx dt = \int \int
\left ( \frac{|\widehat{u_1}|}{\langle \tau - \xi^3 \rangle ^{b + \epsilon}}
\right ) \left ( \langle \tau - \xi ^3 \rangle ^{b + \epsilon} |\widehat{u_2}|
\right ) d\xi d\tau \leq \| u_1 \|_{X^{0,-b-\epsilon}}
\| u_2 \|_{X^{0,b+\epsilon}}
\label{C-S}
\end{equation}

\item The KdV bilinear estimate \cite{KPV}: 
\begin{equation}\| \p_x (u_1 u_2) \|_{X^{0, b}}  \leq \|u_1\|_{X^{0, b + 1}}
\|u_2\|_{X^{0, b + 1}}, \label{bilin}
\end{equation} 
for $b = -\frac12 + \epsilon$, for any $\epsilon > 0$.
     
\item The Strichartz Estimate \cite{KPV2}:
\begin{equation}\|D^{\frac{\theta \alpha}2}_x u \|_{L^q_t L^p_x}  \leq C(\theta,
      \alpha) \|u \|_{X^{0, \frac12 + \epsilon}}, \label{strichartz}
\end{equation} 
for all $(\theta, \alpha) \in [0,1]\times[0,\frac12] $, with 
$p = \frac2{1-\theta}$ and $q = \frac6{\theta(\alpha
         +1)}.$  Here, as elsewhere,  the norm $\| \cdot \|_{L^q_t
L^p_x}$ will mean to
take the $L^p$ norm with respect to $x$ first and then to take the
$L^q$ norm with respect to $t$.
\end{enumerate}

\section {A First Pass at the Theorem}

In this section, we obtain a weaker version of the main result of this
paper.  As mentioned in the introduction, we will follow the structure
of \cite{CKSTT3} because we believe that in this way the argument can
be better understood.  Even though the structure is the same we have
to repeat most of the arguments because the estimates are different.
In addition, in following sections we will use the estimates proven
here.
 
\begin{prop}\label{prop1} Let $ 0 \leq s <1$.  Let $\sigma = \mathrm{dist}_{H^s} (u_0, \Sigma) \ll 1$, and let $u$ be the solution to the KdV such that $u ( \cdot, 0) = u_0$.  Then
$$\mathrm{dist}_{H^s} (u(t), \Sigma) \leq C t^{\frac{1-s}{3-2s-\epsilon}}\sigma ^\frac1{3-2s-\epsilon}$$
for some small $\epsilon > 0$ and for all $t$ such that $t \ll \sigma^{-\frac1{(1-s)-\epsilon}}$.
\end{prop}

\begin{Proof}
 Fix $s$, $u_0$, and $\sigma$.

Let $N \gg 1$. We will fix $N$ later subject to some future
constraints.  Let $I_N$ be the multiplier operator discussed in \ref{I}
with cutoff point $N$.  From now on we will refer to $I_N$ simply as
$I$ unless that is unclear. 

Define 
\begin{equation}
E_N(t) = \Ly (I u (t)) . \label{E}
\end{equation}
Let $\psi$ be a ground state such that $\|u_0 - \psi\|_{H^s} = \sigma$.  
Then $\|Iu_0 - I\psi \| \leq C N^{1-s} \sigma$.
%
Moreover, because $\psi$ is smooth, its Fourier transform is rapidly
decreasing, so $\| I\psi - \psi \|_{H^1} \leq C N^{-C_1}$ for any $C_1$ we
choose.  So, if we require that $N \geq \sigma^{- \epsilon}$ for some
$\epsilon > 0$, then we obtain $\|I\psi - \psi\|_{H^1} \leq C N^{1-s}
\sigma$, so 
$$\| Iu_0 - \psi \|_{H^1} \leq C N^{1-s} \sigma .$$
By (\ref{H1}) this implies (for $\sigma$ sufficiently small with respect
to $N$) that $$|E_N(0) - \Ly(\psi)| \lesssim N^{2 - 2s} \sigma ^2.$$

We will need the following lemma, the proof of which is postponed
until later:

\begin{lemma} \label{lemma1} If there is $t_0 \in \R$ such that $|E_N
(t_0) - \Ly (\psi)| \ll 1$, then $$|E_N (t_0 + \delta) - E_N (t_0)| \leq {\cal{O}}
(\frac1{N^{1 - \epsilon}})$$ where $\delta$ is an absolute constant
depending only on $s$.\footnote{Note that $\delta$ may depend also on
$\|u(t_0)\|_{H^s}$, but $\|u(t_0) - \psi \|_{H^s} \ll 1$ by \cite{W}
and $\| \psi \|_{H^s}$ depends only on $s$, so $\|u(t_0)\|_{H^s}$ also
can be controlled by a constant dependent only on $s$.}
\end{lemma}
For now we will assume the lemma.  Once we have this lemma, by the
same argument which appears in \cite{CKSTT3}, we can iterate to say
that 
$$|E_N(t) - \Ly(\psi)|\lesssim N^{2 - 2s} \sigma ^2,$$ for all $t$ such
that $t \ll N^{1 - \epsilon} N^{2 - 2s} \sigma ^2$.  We may therefore
conclude that, for all such $t$, 
$$\|u(t) - \psi \|_{H^s} \lesssim N^{1 - s} \sigma.$$  
We finally optimize $N$ under the necessary constraints:
\begin{align}
N^{3 - 2s - \epsilon} \sigma ^2 \gg 1 \ \ \ \ \ \ \ \ & \ \ \ \ \ \ \
\ t \ll N^{1 -
\epsilon}N^{2 -2s} \sigma ^2 & N^{2 - 2s} \sigma ^2 \ll 1 \label{conditions1}
\end{align}
and conclude that 
$$\|u - \psi \|_{H^s} \lesssim
t^{\frac{1-s}{3-2s-\epsilon}}\sigma^{\frac1{3 - 2s - \epsilon}},$$
for all $t \ll \sigma^{\frac{-1}{1 - s -\epsilon}}$.   
\end{Proof}

It remains for us to prove the lemma:
 
\begin{Proof} (of Lemma~\ref{lemma1})
To prove Lemma~\ref{lemma1}, we first control $\|Iu\|_{X^{1, \frac12 +
\epsilon}_{[t_0 - \delta, t_0 + \delta]}}$.\footnote{Due to the special
features of the KdV equation, the $X^{s,b}$ norms have been found to
be effective to
work with.}  We will then use this control to take a $\delta$-step
forward in time and measure the growth of $E_N(t)$ in this time step.   

\begin{claim}There exists a $\delta > 0$ such that, for $0 < \epsilon
\ll 1$,  $$\|Iu\|_{X^{1, \frac12 + \epsilon}_{[t_0 - \delta, t_0 +
\delta]}} \lesssim 1.$$
\end{claim}
\begin{Proof} (of Claim)
First note that by (\ref{H1}), $\|Iu(t_0)\|_{H^1}
\lesssim 1$, because $|E_N(t_0) - \Ly(\psi)| \ll 1$ and $\|\psi\|_{H^1}$ is a
constant. 
Moreover, $I$ commutes with differentiation and with $W(t)$.  We may
therefore apply 
the standard $X^{s, b}$ estimates (see, e.g., \cite{KPV}, 
pp. 587-8).  Let $\phi (t) $ be a cutoff
function with support in $[t_0-3, t_0 + 3]$, such that $\phi \equiv 1$ inside
$[t_0-2, t_0+2]$.  Then, $u$ is a fixed point of the operator 
$$L u = \phi(t)W(t - t_0) u(t_0) + \phi(t) \int_{t_0}^t W(t - t')\p_x
(u(t')^2) dt'$$  
on the interval $[t_0 - 2, t_0 + 2]$.  
Then, for $0 < \delta \ll 1$:
\begin{align*}
\|Iu\|_{X^{1, \frac12 + \epsilon}_{[t_0 - \delta, t_0 +
\delta]}} & \leq \|I \psi(t) W(t - t_0) u({t_0}) \|_{X^{1, \frac12 +
\epsilon}_{[t_0 - \delta, t_0 + \delta]}} + \| I \psi(t) \int_{t_0}^t
W(t - t')\p_x (u(t')^2) 
dt'\|_{X^{1, \frac12 + \epsilon}_{[t_0 - \delta, t_0 + \delta]}} \\
& \hskip-.3in \leq C \|Iu(t_0)\|_{H^1} + C \|\p_xIu^2\|_{X^{1, -\frac12 +
\epsilon}_{[t_0 - \delta, t_0 + \delta]}} \\
& \hskip-.3in \leq C \|Iu(t_0)\|_{H^1} + C \delta ^\epsilon \|\p_xIu^2\|_{X^{1, -\frac12 +
2 \epsilon}_{[t_0 - \delta, t_0 + \delta]}}.
\end{align*}
Now, by the bilinear estimate for the KdV (see \cite{KPV}), we have,
for $s \geq -\frac34$:
$$\|\p_x u^2\|_{X^{s, -\frac12 + 2\epsilon}_{[t_0 - \delta , t_0 +
\delta]}} \leq \|u\|^2_ {X^{s, \frac12 + \epsilon}_{[t_0 - \delta, t_0 +
\delta]}}. $$
Consider the multiplier operator $I_1$, which is the same type of
operator as $I$ but with $N = 1$.  It is clear that $\|f\|_{X^{s, b}} \sim
\|I_1 f \|_{X^{1, b}}$, so 
$$ \|\p _x I_1 u^2 \|_{X^{1, -\frac12 +
2 \epsilon}_{[t_0 - \delta, t_0 + \delta]}}  \sim \|\p_x u^2 \|_{X^{s,
-\frac12 + 2 \epsilon}_{[t_0 - \delta, t_0 + \delta]}} \lesssim
\|u\|^2_{X^{s, \frac12 + \epsilon}_{[t_0 - \delta, t_0 + \delta]}}
\sim \|I_1 u \|^2_{X^{1, \frac12 + \epsilon}_{[t_0 - \delta, t_0 + \delta]}}.$$
 But then, by Lemma 12.1 of \cite{CKSTT2} it also follows for general
$N$ that $$\|\p _x I_N u^2 \|_{X^{1, -\frac12 +
2 \epsilon}_{[t_0 - \delta, t_0 + \delta]}}  \lesssim \|I_N u
\|^2_{X^{1, \frac12 + \epsilon}_{[t_0 - \delta, t_0 + \delta]}}.$$
We may therefore conclude that 
$$
\|Iu\|_{X^{1, \frac12 +
\epsilon}_{[t_0 - \delta, t_0 + \delta]}} \leq C \|Iu_0\|_{H^1} + C
\delta ^\epsilon \|Iu\|^2_{X^{1, \frac12 +
\epsilon}_{[t_0 - \delta, t_0 + \delta]}} \leq C + C \delta ^\epsilon \|Iu\|^2_{X^{1, \frac12 +
\epsilon}_{[t_0 - \delta, t_0 + \delta]}} .
$$
Therefore, by a continuity argument, there exists a $\delta > 0$ for which 
$$\|Iu\|_{X^{1, \frac12 +
\epsilon}_{[t_0 - \delta, t_0 + \delta]}} \leq 2C \lesssim
\|Iu(t_0)\|_{H^1} \lesssim 1.$$
This concludes the proof of the claim. \end{Proof}

We now want to take a step forward in time.  Let $f \in H^1$.  Define
$\Omega(f(t)) = \p_t(\Ly(f(t)))$.  Then:
\begin{align}
\Omega(f(t)) = \p_t (\Ly(f)) & = \p_t \left ( \int_\R (f_x^2 + f^2
-\frac23f^3)dx \right ) \notag \\
& = 2 \int_\R (f_x f_{xt} + f f_t - f^2 f_t) dx \notag \\
& = 2 \int_\R f_{t} (-f_{xx} + f - f^2) dx . \label{Omega}
\end{align}
In our case, we are interested in 
\begin{align*}
E_N(t_0 + \delta) - E_N(t_0) & = \int_{t_0}^{t_0 + \delta}
\Omega(Iu(t)) \ dt \\ 
& = 2 \int_{t_0}^{t_0 + \delta} \int_\R Iu_{t} (-Iu_{xx} + Iu - (Iu)^2)
\ dx \ dt \\
& = - 2 \int_{t_0}^{t_0 + \delta} \int_\R (I u_{xxx} +
I(u^2)_x)(-Iu_{xx}+ Iu - (Iu)^2) \ dx \ dt \\
& = - 2 \int_{t_0}^{t_0 + \delta} \int_\R Iu_{xxx}((Iu)^2 - Iu^2) \ dx
\ dt + 2 \int_{t_0}^{t_0 + \delta} \int_\R (u^2)_x I^2 u \ dx \ dt \ + \\
& \phantom{moveit!} - 2
\int_{t_0}^{t_0 + \delta} \int_\R I(u^2)_x (Iu)^2 \ dx \ dt \\
& = -2 \int_{t_0}^{t_0 + \delta} \Lambda_3\left (\xi_1^3
m(\xi_1)(m(\xi_2)m(\xi_3) - m(\xi_2 + \xi_3) ) ; u ; u ; u
\right ) \ dt \ + \\
& \phantom{moveit!} + 4 \int_{t_0}^{t_0 + \delta} \Lambda_3\left ( \xi_1 m(\xi_3)^2 ) ; u ; u ; u \right ) \ dt \ + \\
& \phantom{moveit!} - 4 \int_{t_0}^{t_0 + \delta}
\Lambda_4\left ( \xi_1 m(\xi_1 + \xi_2) m(\xi_3)m(\xi_4) ; u
; u ; u ; u \right ) \ dt .
\end{align*}

We will prove the following more general estimates in order to
control $E_N(t_0~+~\delta)~-~E_N(t_0)$:
\begin{align}
\left | \int_{t_0}^{t_0 + \delta} \Lambda_3\left (\xi_1^3
m(\xi_1)(m(\xi_2)m(\xi_3) - m(\xi_2 + \xi_3) ) ; u_1 ; u_2 ; u_3
\right ) \ dt \right | & \lesssim N^{-1 + \epsilon} \prod_{i = 1}^{3} \|Iu_i \|_{X^{1,
\frac12 + \epsilon}_{[t_0 - \delta, t_0 + \delta]}},  \label {estimate1}\\
\left | \int_{t_0}^{t_0 +\delta} \Lambda_3\left ( \xi_1 m(\xi_3)^2 ) ;
u_1 ; u_2 ; u_3 \right ) \ dt \right |
& \lesssim  N^{-1 + \epsilon} \prod_{i = 1}^{3} \|Iu_i \|_{X^{1,
\frac12 + \epsilon}_{[t_0 - \delta, t_0 + \delta]}},  \label{estimate2}\\
\left |\int_{t_0}^{t_0 +\delta} \Lambda_4\left ( (\xi_1)m(\xi_1 + \xi_2)
m(\xi_3)m(\xi_4) ; u_1 
; u_2 ; u_3 ; u_4 \right ) \ dt \right | & \lesssim N^{-1 + \epsilon}
\prod_{i = 1}^{4} \|Iu_i \|_{X^{1, 
\frac12 + \epsilon}_{[t_0 - \delta, t_0 + \delta]}} . \label{estimate3}
\end{align}
Recall that $m(\xi)$ is the multiplier associated with the operator
$I$, and it is identically $1$ for $|\xi| \leq N$, and equals
$\frac{\xi^{s - 1}}{N^{s-1}}$ for $|\xi| > 10N$.  Note that because
our norms are of $L^2$ type, we may replace 
$\hat{u}$ by $|\hat{u}|$ without affecting the estimates.  
For each estimate, we will divide the functions $u_i$ into dyadic
blocks $N_i$ in frequency space and make appropriate estimates.  We will
then sum over these dyadic blocks to obtain the full estimate.

\medskip
\begin{Proof} (of Estimate (\ref{estimate1})) We consider the multiplier
$N_1^3 m(N_1) (m(\xi_2)m(\xi_3) - m(\xi_2 + \xi_3))$. 
Recall that we have $N_1 + N_2 + N_3 = 0$ and note that we may assume
$N_2 \geq N_3$ because of the symmetry, and that $N_2 > N$ or else the
whole symbol is $0$.  We will consider two cases:
\begin{enumerate}
\item $N_2 \gg N_3$: This implies that $N_1 \sim N_2$.  

First suppose that $N_3 \leq N$.  Then $m(N_3) = 1$, so 
$$N_1^3 m(N_1) (m(N_2)m(N_3) - m(N_2 + N_3)) = N_1^3 m(N_1)
(m(N_2) - m(N_2 + N_3)).$$  By the mean value theorem, this is $\leq
N_1^3 m(N_1) m'(N_2) N_3$, so, since $N_1 \sim N_2$ and $m'(N_2) =
\frac{m(N_2)}{N_2}$, 
\begin{equation}\label{estimate1a}
\N_1^3 m(N_1) (m(N_2)m(N_3) - m(N_2 + N_3)) leq N_1 N_2 N_3 m(N_1)
m(N_2) m(N_3). 
\end{equation} 
Now, consider the whole integral, and use inequality
(\ref{C-S}):\footnote{Here we are ignoring that we are on a finite
time interval.  To be precise one should repeat the argument given in
\cite{CKSTT3b} during the proof of Lemma 8.1}
\begin{align*}
N_1 N_2 N_3 & m(N_1) m(N_2) m(N_3) \int_{t_0}^{t_0 + \delta} \int | \hat{u}_{1,N_1} \hat{u}_{2,
N_2} \hat{u}_{3, N_3} | d\xi \ dt \\
& \leq N_1 N_2 N_3  m(N_1) m(N_2) m(N_3)
\|u_{1,N_1} u_{3, N_3} \|_{X^{0, -\frac12 + \epsilon}}
\|u_{2,N_2}\|_{X^{0, \frac12 + \epsilon}} .
\end{align*}
Then, by the KdV bilinear estimate and because $\| \p_x u_{1, N_1} \| \sim
N_1 \|u_{1, N_1} \|$ and $N_3 \ll N_1$, we obtain:
\begin{align*}
N_1 N_2 N_3 & m(N_1) m(N_2) m(N_3) \int_{t_0}^{t_0 + \delta} \int | \hat{u}_{1,N_1)} \hat{u}_{2,
N_2} \hat{u}_{3, N_3} | d\xi \ dt\hskip1.5in \\
& \leq N_1 N_2 N_3 m(N_1) m(N_2) m(N_3)
\frac1{N_1} \|\p_x (u_{1, N_1} u_{3, N_3} )\|_{X^{0, -\frac12 +
\epsilon}} \|u_2\|_{X^{0, \frac12 + \epsilon}}\\
& \leq N_2 N_3 m(N_1) m(N_2) m(N_3) 
\|u_{1, N_1}\|_{X^{0, \frac12 + \epsilon}} \|u_{2, N_2}\|_{X^{0, \frac12 +
\epsilon}} \|u_{3, N_3}\|_{X^{0, \frac12 + \epsilon}} .
\end{align*}
But then, by definition of $I$ and the $X^{s,b}$ spaces and because
$N_1 \sim N_2$, this is
controlled by 
$$\frac1{N_1^{\frac12 - \tilde{\epsilon}}}
\frac1{N_2^{\frac12}} \frac1{N_3^{\tilde{\epsilon}}} \|Iu_{1, N_1} \|_{X^{1, \frac12
+\epsilon}} \|Iu_{2, N_2} \|_{X^{1, \frac12 +\epsilon}} \|Iu_{3,
N_3}\|_{X^{1, \frac12 +\epsilon}} .$$ 
When we sum this in the $N_i$s, we will lose a power of
$\epsilon$, and obtain a term of size $\frac1{N^{1 - \epsilon}}$ as
claimed. 

\medskip

Now, suppose instead that $N_2 \gg N_3 > N$.  Then 
\begin{align*} 
N_1^3 m(N_1) & (m(N_2) m(N_3) - m(N_2 + N_3)) = \\
& N_1^3m(N_1)(m(N_2)m(N_3) - m(N_3) m(N_2 + N_3)) + N_1^3m(N_1)(m(N_3)
m(N_2 + N_3) - m(N_2 + N_3)) \\
& = M_1 + M_2 
\end{align*}
For estimate $M_1$, use the mean value theorem (recall that $N_1 \sim N_2$):
$$ M_1 \leq N_1 N_2 N_3 m(N_1)m(N_2) m(N_3).$$
Then the same calculation as before implies
that the part of the left-hand side of (\ref{estimate1}) containing
$M_1$ also sums to $\frac1{N^{1 - \epsilon}}$ as desired. 

On the other hand, 
$$M_2 = N_1^3m(N_1)m(N_2 + N_3)(m(N_3) - 1).$$
Note that $|m(N_3) - 1| \leq 2$, and $m(N_2  + N_3) \sim m(N_2)$ because $N_2
\gg N_3$.  So, 
$$ M_2 \leq N_1^3 m(N_1) m(N_2) \frac{m(N_3)}{m(N_3)} \lesssim
\frac{N_1^2 N_2 N_3 
m(N_1) m(N_2) m(N_3)}{N_3 m(N_3)}, $$
where $m(N_3) \sim \frac{N_3^{s-1}}{N^{2-1}}$, so $\frac1{N_3 m(N_3)}
\sim \frac1{N_3^s N^{1-s}}$.  Therefore, we find that 
 
$$M_2 \lesssim \frac1{N^{1-s}}\frac1{N_3^{s}}N_1^2 N_2 N_3 m(N_1) m(N_2)
m(N_3) N_1.$$
As before, we compute that
\begin{align*}
\int_{t_0}^{t_0 + \delta} \int | \hat{u}_{1,N_1} \hat{u}_{2,N_2}
\hat{u}_{3, N_3} | d\xi \ dt & \leq 
\|u_1 u_3 \|_{X^{0, -\frac12 + \epsilon}} \|u_2\|_{X^{0, \frac12 +
\epsilon}} \\
& \leq \frac1{N_1} \|u_{1, N_1}\|_{X^{0, \frac12 + \epsilon}} \|u_{2, N_2}\|_{X^{0, \frac12 +
\epsilon}} 
\|u_{3, N_3}\|_{X^{0, \frac12 + \epsilon}}. 
\end{align*}
And so, the part of the left-hand side of (\ref{estimate1}) containing
$M_2$ is bounded by
$$\frac1{N^{1-s}}\frac1{N_3^{s}} \|Iu_{1, N_1} \|_{X^{1, \frac12
+\epsilon}} \|Iu_{2, N_2} \|_{X^{1, \frac12 +\epsilon}} \|Iu_{3,
N_3} \|_{X^{1, \frac12 +\epsilon}}.$$
To sum this, we use Cauchy-Schwartz and the fact that $N_1 \sim N_2$
, to obtain the same estimate as before.

\item Now consider the case where $N_2 \sim N_3$.  Then $N_1 = -(N_2 +
N_3)$ may be smaller.  We once again want to estimate the multiplier 
$$N_1^3 m(N_1) (m(N_2)m(N_3) - m(N_2 + N_3)) = N_1^3 m(N_1) m(N_2)
m(N_3) - N_1^3 m(N_1) ^2 = M_3 + M_4.$$
We have 
\begin{align*}
M_3 = N_1^3 m(N_1)m(N_2)m(N_3) & = N_1 N_2 N_3 m(N_1)m(N_2)m(N_3)
\frac{N_1^2}{N_2 N_3} \\
& \lesssim  N_1 N_2 N_3 m(N_1)m(N_2)m(N_3).
\end{align*}
Then, by the same argument as for the first part of the first case,
this sums to ${\cal{O}} (\frac1{N^{1-\epsilon}})$.
For $M_4$, we have 
\begin{align*}
M_4 & \lesssim N_1 N_2 N_3 m(N_1)m(N_2)m(N_3)
\frac{N_1^2}{N_2 N_3} \frac{m(N_1)}{m(N_2)m(N_3)} \\
& = N_1 N_2 N_3
m(N_1)m(N_2)m(N_3) \frac{N_1^{1 - s}}{N^{1 - s} N_2^s N_3^s}.
\end{align*}
We then use the bilinear estimate as before to conclude that 
\begin{align*}
N_1^3 m (N_1)^2 & \int_{t_0}^{t_0 + \delta} \int | \hat{u}_{1,N_1}
\hat{u}_{2,N_2} \hat{u}_{3, N_3} | d\xi \ dt \\
& \lesssim N_1 N_2 N_3 m(N_1) m(N_2) m(N_3)
\frac{N_1^s}{N^{1-s}N_2^s N_3^s} \|u_{1, N_1}\|_{X^{0, \frac12 +
\epsilon}} \|u_{2, N_2}\|_{X^{0, \frac12 + \epsilon}} \|u_{3,
N_3}\|_{X^{0, \frac12 + \epsilon}} \\
& \lesssim \frac1{N^{1-s}} \frac1{N_1^{2
\epsilon}} \frac1{N_2^{\frac{s}2 - \epsilon}} \frac1{N_3^{\frac{s}2 -
\epsilon}}  
\|Iu_{1, N_1} \|_{X^{1, \frac12
+\epsilon}} \|Iu_{2, N_2} \|_{X^{1, \frac12 + \epsilon}} \|Iu_{3,
N_3} \|_{X^{1, \frac12 +\epsilon}}, 
\end{align*}
after using again the fact that $N_1 \leq N_2 \sim N_3$.  Summing in
the $N_i$s, we can see that this again gives ${\cal{O}} (\frac1{N^{1-\epsilon}})$.
\end{enumerate}
This concludes the proof of estimate (\ref{estimate1}).
\end{Proof}

We next need to prove the estimate (\ref{estimate2}):
$$\left |\int_{t_0}^{t_0 + \delta} \Lambda_3\left ( \xi_1 m(\xi_3)^2 ) ; u ;
u ; u \right ) \ dt \right |
\lesssim  N^{-1 + \epsilon} \prod_{i = 1}^{3} \|Iu_i \|_{X^{1,
\frac12 + \epsilon}} 
$$
\begin{Proof} (of Estimate (\ref{estimate2}))

We will consider the multiplier $N_1 m(N_3)^2$.  Note that if $N_1$,
$N_2$, and $N_3$ are all less than $N$, then the 
operator given by the symmetrization of this multiplier is identically
zero.  So at least one of $N_1$, $N_2$, and $N_3$ must be greater than
$N$. If $N_3 <
N$, this multiplier is just $N_1$, and, as above, 
\begin{align*}
N_1 \int_{t_0}^{t_0 + \delta} \int | \hat{u}_{1,N_1} \hat{u}_{2,N_2}
\hat{u}_{3, N_3} | d\xi \ dt & \lesssim
N_1 \|u_{1,N_1} u_{3,N_3} \|_{X^{0, -\frac12 + \epsilon}} \|u_{2,N_2}\|_{X^{0, \frac12 +
\epsilon}} \\
& \lesssim \|u_{1, N_1}\|_{X^{0, \frac12 + \epsilon}} \|u_{2,
N_2}\|_{X^{0, \frac12 + \epsilon}} \|u_{3, N_3}\|_{X^{0, \frac12 +
\epsilon}} \\
& \lesssim \frac1{N_1 m(N_1) N_2 m(N_2) N_3 m(N_3)} \prod_{i = 1}^3 \|Iu_{i,
N_i} \|_{X^{1, \frac12 + \epsilon}} .
\end{align*}
Since at least one
of $N_1$, $N_2$ is greater than $N$, the quantity computed above sums to
no more than ${\cal O} (\frac1{N^{1-\epsilon}})$.

Now, if $N_3 > N$, as above
$$ N_1 m(N_3)^3 \int_{t_0}^{t_0 + \delta} \int | \hat{u}_{1,N_1}
\hat{u}_{2,N_2} \hat{u}_{3, N_3} | d\xi \ dt \lesssim \frac{m(N_3)}{N_1 m(N_1) N_2 m(N_2) N_3} \prod_{i = 1}^3 \|Iu_{i,
N_i} \|_{X^{1, \frac12 + \epsilon}} ,
$$ which again sums to ${\cal O} (\frac1{N^{1-\epsilon}})$ in
the worst cases.  
\end{Proof}

Finally, we need to prove estimate (\ref{estimate3}):
$$\left | \int_{t_0}^{t_0 + \delta} \Lambda_4\left ( (\xi_1)m(\xi_1 +
\xi_2) m(\xi_3)m(\xi_4) ; u 
; u ; u ; u \right ) \ dt \right | \lesssim N^{-1 + \epsilon} \prod_{i =
1}^{4} \|Iu_i \|_{X^{1, 
\frac12 + \epsilon}} $$
\begin{Proof} (of Estimate (\ref{estimate3}))
We consider the multiplier $N_1 m(N_1 + N_2)m(N_3) m(N_4)$.  Recall
that we have $N_1 + N_2 + N_3 + N_4 = 0$ and by symmetry we may
assume $N_3 \geq N_4$.  
Consider $$\int | \hat{u}_{1, N_1}   \hat{u}_{2, N_2} \hat{u}_{3,
N_3}\hat{u}_{4, N_4} | d\xi \lesssim \prod_{i = 1}^4 \|u_{i, N_i}
\|_{L^4} . $$
We use the Strichartz estimate (\ref{strichartz}) with $(\theta,
\alpha) = (\frac12, 0)$ and $p = 4$, $q = 12$, obtaining
$$\|u\|_{L^{12}_tL^4_x} \leq C \|u\|_{X^{0, \frac12 + \epsilon}}.$$
In our case, therefore, we may conclude that
\begin{align*}
\int_{t_0}^{t_0 + \delta} \int | \hat{u}_{1, N_1}   \hat{u}_{2, N_2} \hat{u}_{3,
N_3}\hat{u}_{4, N_4} | d\xi dt & \leq \int_{t_0}^{t_0 + \delta} \prod_{i = 1}^4 \|u_{i, N_i}
\|_{L^4} \ dt \\
& \leq \|1\|_{L^\frac32_t} \prod_{i = 1}^4 \|u_{i, N_i}\|_{L^{12}_t
L^4_x} \ \leq \ C \delta^{\frac23} \prod_{i = 1}^4 \|u_{i, N_i}\|_{X^{0,
\frac12 + \epsilon}} .
\end{align*}
Therefore 
\begin{align*}
N_1 m(N_1& + N_2) m(N_3) m(N_4) \int | \hat{u}_{1, N_1}   \hat{u}_{2, N_2} \hat{u}_{3,
N_3}\hat{u}_{4, N_4} | d\xi \hskip1.5in \\
& \leq N_1 N_2 N_3 N_4 m(N_1) m(N_2) m(N_3) m(N_4)
\frac{m(N_1 + N_2)}{N_2 N_3 N_4 m(N_1) m(N_2)} \prod_{i = 1}^4 \|u_{i,
N_i}\|_{X^{0, \frac12 + \epsilon}} \\
& \leq \frac{m(N_1 + N_2)}{m(N_1) m(N_2)}\frac1{N_2 N_3 N_4}
\prod_{i = 1}^4 \|Iu_{i, N_i}\|_{X^{1, \frac12 + \epsilon}} .
\end{align*}
We will now estimate $\frac{m(N_1 + N_2)}{m(N_1) m(N_2)}\frac1{N_2 N_3
N_4}$, considering several cases (recall that $N_3 \geq N_4$): 
\begin{enumerate}
\item First assume $N_1 \gg N_2$.  

If $N_1 \leq N$, $m(N_1 + N_2) =
m(N_1) = m(N_2) = 1$.  Note that if $N_1$, $N_2$, $N_3$, $N_4$ are all
less than $N$, then the operator is identically zero by
symmetrization.  Hence at least one of the dyadic blocks must be at
least $N$ for the operator to be nontrivial.  Therefore, if $N_1 \leq
N$, then $N_3 > N$.  Hence the multiplier, which reduces to
$\frac1{N_2 N_3 N_4}$ in this case clearly sums to no more than
${\cal O}(\frac1{N^{1-\epsilon}})$.  

So we may assume that $N_1 > N$.  Then $m(N_1 + N_2) \sim m(N_1)$
because $N_1 \gg N_2$.  So our multiplier reduces to $\frac1{m(N_2)}
\frac1{N_2 N_3 N_4}$.  If $N_2 < N$, this is again $\frac1{N_2 N_3
N_4}$.  But now, because $N_1 + N_2 + N_3 + N_4 = 0$, $N_1 \sim N_3$.
Hence we may write
$$\frac1{N_2 N_3 N_4} \leq \frac1{N_1^{\frac12} N_2 N_3^{\frac12}
N_4}, $$
which sums to ${\cal O} (\frac1{N^{2 - \epsilon}})$.  

Finally, if $N_2 > N$ as well, we have $m(N_2) \sim \frac{N_2^{s-1}}{N^{s-1}}$.
Therefore, because $N_1 \sim N_3$, the multiplier is controlled by 
$$\frac1{N^{1-s}N_2^s N_1^\frac12 N_3^\frac12 N_4}, $$
which sums to ${\cal O} (\frac1{N^{3 - \epsilon}})$.

\item $N_1 \ll N_2$.

Then $m(N_1 + N_2) \sim m(N_2)$, so the multiplier is 
$$\frac1{m(N_1) N_2 N_3 N_4}.$$
The case where $N_1$ and $N_2$ are both less than $N$ is the same as
before.  So we consider first what happens when $N_1 < N$.  Then we
again have $\frac1{N_2 N_3 N_4}$.  As before, the operator is trivial
unless $N_3 > N$, and when $N_3 > N$ this sums to ${\cal O}
(\frac1{N^{1 - \epsilon}})$ as desired.  

If instead $N_1 > N$, we have $$\frac{N_1^{1 - s}}{N^{1 - s} N_2 N_3
N_4} \leq \frac{1}{N^{1 - s} N_1^\epsilon N_2^{\frac{1 + s -
\epsilon}2} N_3^{\frac{1 + s - \epsilon}2} N_4} , 
$$ which sums to ${\cal O} (\frac1{N^{3 - \epsilon}})$ as in the first
case.

\item Finally we consider the case where $N_1 \sim N_2$.  

Once again the case where both $N_1$ and $N_2$ are less than $N$ is
the same as before.  Therefore, we consider the case where $N_1 \sim
N_2 > N$.  Then $m(N_1 + N_2) \leq 1$, so the multiplier reduces to $\frac1{m(N_1) m(N_2) N_2 N_3 N_4}.$  If $N_3 \sim N_1
\sim N_2$, then this is controlled by $\frac1{N^{2(1-s)} N_2^s N_3^{s}
N_4}$ and since $N_3$ controls all the other quantities, we may again
sum to conclude that this is bounded by ${\cal O} (\frac1{N^{2 - 
\epsilon}})$.  

We must at last consider the case $N_3 \ll N_1$.  For this case we
must reconsider the original calculations done at the beginning of
this estimate.  Instead of treating all four functions equally, we
will write:
$$\int | \hat{u}_{1, N_1}   \hat{u}_{2, N_2} \hat{u}_{3,
N_3}\hat{u}_{4, N_4} | d\xi \leq \|u_{1,N_1} u_{3,N_3}\|_{L_2} \|u_{2, N_2}
\|_{L^4}\|u_{4, N_4} \|_{L^4} . $$
Therefore, using the Strichartz estimate again, the fact
that $\|f\|_{X^{0,0}} \leq \|f\|_{X^{0, \frac12 + \epsilon}}$ for any
function $f$, and the KdV bilinear estimate:
\begin{align*}
\int_{t_0}^{t_0 + \delta} \int | \hat{u}_{1, N_1}   \hat{u}_{2, N_2}
\hat{u}_{3, 
N_3}\hat{u}_{4, N_4} | d\xi dt & \leq \int_{t_0}^{t_0 + \delta}
\|u_{1,N_1} u_{3,N_3}\|_{L^2} \|u_{2, N_2} \|_{L^4}\|u_{4, N_4}
\|_{L^4} \ dt \\ 
& \leq \|1\|_{L^3_t} \|u_{1,N_1} u_{3,N_3}\|_{L^2_tL^2_x} \|u_{2, N_2}
\|_{L^{12}_tL^4_x}\|u_{4, N_4} \|_{L^{12}_tL^4_x} \\
& \leq \ C \delta^{\frac13} \|u_{1,N_1} u_{3,N_3}\|_{X^{0, \frac12
+\epsilon}} \|u_{2, N_2}\|_{X^{0,
\frac12 + \epsilon}} \|u_{4, N_4}\|_{X^{0, \frac12 + \epsilon}} \\
& \leq C \delta^\frac13 \frac1{N_1 + N_3} \prod_{i=1}^4 \|u_{i, N_i}\|_{X^{0,
\frac12 + \epsilon}} .
\end{align*}
Now, recall that $N_3 \ll N_1$, $N_3 \geq N_4$, and $N_1 + N_2 + N_3 +
N_4 = 0 $, so $N_3 + N_1 \sim N_2 \sim N_1$.  Therefore, our entire
operator may be estimated as follows:
\begin{align*}
N_1 m(N_1 & + N_2) m(N_3) m(N_4) \int | \hat{u}_{1, N_1}   \hat{u}_{2,
N_2} \hat{u}_{3, N_3} \hat{u}_{4, N_4} | d\xi \hskip1.5in \\
& \leq m(N_1 + N_2) m(N_3) m(N_4) \prod_{i=1}^4 \|u_{i, N_i}\|_{X^{0,
\frac12 + \epsilon}} \\
& \leq \frac1{N_1 N_2 N_3 N_4} \frac{m(N_1 + N_2)}{m(N_1) m(N_2)}
\prod_{i=1}^4 \|I u_{i, N_i}\|_{X^{1, \frac12 + \epsilon}} .
\end{align*}
We therefore need only to sum $$\frac1{N_1 N_2 N_3 N_4} \frac{m(N_1 +
N_2)}{m(N_1) m(N_2)} \lesssim \frac1{N^{2(1-s)} N_1^s N_2^s N_3 N_4}$$
which as before is at worst ${\cal O} (\frac1{N^{2 - \epsilon}})$.
\end{enumerate}
This concludes the proof of estimate (\ref{estimate3}).
\end{Proof}

Having proved all three estimates, we note that 
$$|E_N(t_0 + \delta) - E_N(t_0)| \lesssim 2( (\ref{estimate1}) -
(\ref{estimate2}) + (\ref{estimate3})) \leq {\cal O} (\frac1{N^{1 -
\epsilon}})$$
since we have already checked that $\|Iu\|_{X^{1, \frac12 +
\epsilon}_{t_0 - \delta, t_0 + \delta}} \lesssim 1$.  This concludes
the proof of Lemma~\ref{lemma1}. 
\end{Proof}

\section{A Second Pass at the Theorem}\label{section2}
In this section, we will improve the powers of $t$ and of $\sigma$ which
appear in Proposition~\ref{prop1}.  We will do this by exploiting more
carefully the fact that $\|u_0 - \psi\|_{H^s}$ is small.  

\begin{prop} \label{prop2} Let $0 \leq s < 1$ and suppose
$\mathrm{dist}_{H^s} (u_0, \Sigma) = \sigma \ll 1$.  Then we have, for
some small $\epsilon > 0$,  
$$\mathrm{dist}_{H^s} (u(t) , \Sigma) \leq t^{1 - s + \epsilon}
\sigma^{1 + \epsilon}$$
for all $t$ such that $1 < t \ll \sigma^{-\frac1{1-s} - \epsilon}$.
\end{prop}

\begin{Proof} Fix $s, u_0,$ and $\sigma$.  We retain the
definition of $E_N(t)$ (see \ref{E}), and the set-up of the previous
proposition. 
The main difference will be a sharper estimate for 
$E_N(t_0~+~\delta)~-~E_N(t_0)$:
\begin{lemma} \label{lemma2} If there is a $t_0 \in \R$ and $\tilde{\sigma}$
with $N^{-C} < \tilde{\sigma} \ll 1$ for some arbitrary constant $C$,
such that for some solution to (\ref{ODE}) $\psi$, $|E_N(t_0) - \Ly (\psi)| \leq \tilde{\sigma} ^2$ then we have,
for some $\delta > 0$ depending only on $s$, 
$$E_N(t_0 + \delta) = E_N(t_0) + {\cal O} (\frac1{N^{1 - \epsilon}}
\tilde{\sigma}^2).$$
\end{lemma}
We will, as in the previous section, postpone the proof of the lemma
until later.  First we will complete the proof of
Proposition~\ref{prop2} taking advantage of Lemma~\ref{lemma2}.  We
can again iterate the lemma.  Let $\tilde{\sigma} = N^{1 - s} \sigma$.  We
then obtain
$$|E_N(t) - \Ly (\psi)| \lesssim N^{2 - 2s} \sigma ^2 , $$
for $1 \leq t \ll N^{1 - \epsilon}$, and by \ref{H1} we can then conclude that
for all such times $t$, $\mathrm{dist}_{H^s} (u(t), \Sigma) \lesssim~N^{1
- s} \sigma$.  But now, we may optimize $N$ under the conditions 
\begin{align}
N^{-C} > \sigma \ \ \ \ \ \ \ \ \ \ & \ \ \ \ \ \ \ \ \ 
\ t \ll N^{1 - \epsilon} & N^{2 - 2s} \sigma ^2 \ll 1 . \label{conditions2}
\end{align}
Contrast these conditions with (\ref{conditions1}).  With this
improvement, we obtain 
$$\mathrm{dist}_{H^s} (u(t), \Sigma) \lesssim t^{1 - s + \epsilon}
\sigma^{1 + \epsilon} , $$ 
for $1 \leq t \ll \sigma^{-\frac1{1 - s + \epsilon}}$ as claimed.
\end{Proof}

It therefore remains only to prove the lemma:

\begin{Proof} (of Lemma~\ref{lemma2})

By \ref{H1} and the calculations at the start of Lemma~\ref{lemma1},
there exists a $\psi \in \Sigma$ such that \linebreak $\|Iu(t_0) - \psi\|_{H^1}
\lesssim \tilde{\sigma}$.  Let $Q(x,t) = \psi(x - t)$.  Define 
$$w(x,t) = u(x,t) - Q(x,t) . $$
As before $\psi$ is Schwartz and since $N^{-C} \lesssim \tilde{\sigma}$
for some $C$, we may conclude that 
$\|Iu(t_0) -~I\psi\|_{H^1} \lesssim~\tilde{\sigma}$, i.e. $\|w(t_0)\|_{H^1} \lesssim \tilde{\sigma}$. 

\begin{claim} $\|Iw\|_{X^{1, \frac12 + \epsilon}_{[t_0 - \delta, t_0 +
\delta]}} \lesssim \tilde{\sigma}.$ \label{claim2}
\end{claim}
\begin{Proof} 
The function $w(t)$ obeys the following difference equation:
\begin{equation}
w_t + w_{xxx} + \p_x(w (w + 2 Q) ) = 0 . \label{diffeqn}
\end{equation}
We can therefore use the standard $X^{s,b}$ estimates as in
Lemma~\ref{lemma1} to conclude that  
$$\|Iw\|_{X^{1, \frac12 + \epsilon}_{[t_0 - \delta, t_0 + \delta]}}
\leq \|Iw(t_0)\|_{H^1} + \delta^{\epsilon} \|I(\p_x (w (w
+2Q))\|_{X^{1, -\frac12 + \epsilon}_{[t_0 - \delta, t_0 + \delta]}} . 
$$ 
We then use the bilinear estimate as in Lemma~\ref{lemma1}, as well as
the fact that  $Q$ is a Schwartz function in $x$, to conclude that 
$$\|Iw\|_{X^{1, \frac12 + \epsilon}_{[t_0 - \delta, t_0 + \delta]}}
\leq \tilde{\sigma} + C \delta ^\epsilon \|Iw\|_{X^{1, \frac12 +
\epsilon}_{[t_0 - \delta, t_0 + \delta]}} + \delta ^\epsilon
\|Iw\|^2_{X^{1, \frac12 + \epsilon}_{[t_0 - \delta, t_0 + \delta]}}
$$
and therefore, by a continuity argument again,
$\|Iw\|_{X^{1, \frac12 + \epsilon}_{[t_0 - \delta, t_0 + \delta]}}
\lesssim \tilde{\sigma}$ for some $\delta > 0$ sufficiently small.
This concludes the proof of the claim.  \end{Proof}

Finally, we must again take a $\delta$-step forward in $t$.  We will
show that 
$$E_N(t_0 + \delta) - E_N(t_0) = 2 \int_{t_0}^{t_0 + \delta}
\Omega(I(Q+ w)(t))
dt = {\cal O} (\frac1{N^{1 - \epsilon}} \tilde{\sigma}^2) . 
$$
We will use Lemma~\ref{lemma1} to do this, following the method of
\cite{CKSTT3}, rather than checking it directly.

Because $\tilde{\sigma} \gtrsim N^{-C}$ it will suffice to prove the
more general bound 
$$E_N(t_0 + \delta) - E_N(t_0)  = {\cal O} (\frac1{N^{1 - \epsilon}}
\tilde{\sigma}^2) + {\cal O} (\frac1{N^{C + 1}} \tilde{\sigma}) +
{\cal O}(\frac1{N^{2C + 1}}) .
$$
To do so, consider $\Omega(I(Q(t) + \frac{k}{\tilde{\sigma}} w(t)))$
for $|k| \leq 1$.
Recall that if $f$ is a solution to the KdV, then
$$\Omega(If(t)) = \langle If_{xxx}, (If)^2 - If^2 \rangle - \langle
(f^2)_x, I^2 f \rangle + \langle I (f^2)_x, (If)^2 \rangle ,
$$ where $\langle \ , \ \rangle$ denotes the $L^2$ inner product.
Therefore, $\Omega(IQ(t) + \frac{k}{\tilde{\sigma}} Iw(t))$ is a
polynomial in $k$.  In addition,
from the estimates in Lemma~\ref{lemma1}, which applies to $I(Q(t) +
\frac{k}{\tilde{\sigma}} w(t))$ 
because $\|\frac{k}{\tilde{\sigma}} Iw(t) \|_{X^{1, \frac12 +
\epsilon}_{[t_0 - \delta, t_0 + \delta]}} \ll 1$ for $|k| < 1$, we may
conclude that the coefficients of the polynomial
$$P_\delta(k) = 2 \int_{t_0}^{t_0 + \delta} \Omega(IQ(t) +
\frac{k}{\tilde{\sigma}} Iw(t)) \ dt$$
are ${\cal O} (\frac1{N^{1 - \epsilon}}) $ so all the terms of second
order or higher will validate the desired inequality automatically.

We therefore need only to check that the constant term is ${\cal O}
(\frac1{N^{2 C + 1}})$ and the linear terms are
${\cal O} (\frac1{N^{C + 1}})$.  
The constant term comes from 
$$\Omega(I(Q(t)) = \langle IQ_t, IQ - IQ_{xx} - (IQ)^2 \rangle =
\langle IQ_t, IQ^2 - (IQ)^2 \rangle,$$
which is ${\cal O}{\frac1{N^{2C + 1}}}$ because $IQ^2 - (IQ)^2 = IQ(I -
1)Q + Q(I - 1)Q + (1 - I)Q^2$. But now note that because $Q$ is
Schwartz and $m(\xi) \equiv 1$ for $|\xi| \leq N$, $m(\xi) \leq 1$ for all $\xi$, we may conclude that $(I - 1) Q = {\cal O} (N^{-2C})$ for
any $C$ we choose because $Q(t)$ is Schwartz in $x$.  The same is true
for $Q^2$.   

For the linear term, note that the linear term of $E_N(t)$ is given by:
$$E_N(t) = 2 \langle IQ(t)_x, Iw(t)_x \rangle + 2 \langle IQ(t), Iw(t)
\rangle - 2 \langle (IQ)^2(t), w(t) \rangle,$$
so the linear term of $\Omega(t)$ is: 
\begin{align*}
\Omega(t) & = \frac{d}{dt} E_N(t) \\
& = 2 \langle w_t, I (IQ^2 - (IQ)^2)
\rangle + 2 \langle w, \frac{d}{dt}\left ( I(IQ^2 - (IQ)^2) \right )
\rangle  + \mbox{higher order terms}\\
& = -2 \langle w_{xxx} , I (IQ^2 - (IQ)^2) \rangle + 2 \langle w,
\frac{d}{dt} (I (IQ^2 - (IQ)^2)) \rangle + \mbox{higher order terms} .
\end{align*}
We can thus bound those linear terms by: (after integrating by parts)
$$\|w\|_{L^2} (\| I(IQ^2 - (IQ)^2) \|_{H^3} + \|\p_t (I(IQ^2 -
(IQ)^2)) \|_{L^2})$$
But now note again that because $Q$ is Schwartz we may conclude that $\|IQ^2 -
(IQ)^2\|_{H^s} \lesssim N^{-C}$ for any $s > 0$.  The same is true for
$Q_t$.  Therefore, the linear terms of $\Omega(t)$ are controlled by
$\|w\|_{L^2} N^{-C-1}$, and so, also using the fact that $\|w\|_{X^{1,
\frac12 + \epsilon}_{[t_0 - \delta, t_0 + \delta]}} \leq \sigma$, we
conclude that 
$$|E_N(t_0 + \delta) - E_N(t_0)| \lesssim {\cal O} (N^{-2C - 1}) + {\cal O} (N^{-C-1}
\tilde{\sigma}) + {\cal O } (N^{-1 + \epsilon} \tilde{\sigma}^2)$$
which concludes the proof of Lemma~\ref{lemma2}.
\end{Proof}

\section{Final Proof of the Main Theorem}

In this section we will at last obtain the full power of
Theorem~\ref{main}:
\begin{theorem}\label{main2}

Let $0 \leq s < 1$, Let $\sigma = \mathrm{dist}_{H^s} (u_0, \Sigma) \ll 1$, and
let $u$ be the solution to the KdV such that $u ( \cdot , 0 ) = u_0$.
Then $\mathrm{dist}_{H^s} (u(t), \Sigma) \leq t^{1-s+\epsilon}
\sigma$, for all $t$ such that $t \ll \sigma^{-\frac1{1-s+\epsilon}}.$

\end{theorem}

To do so, we will need to refine the choice of the soliton $Q$ to which $u$
is close.  In the previous section, we chose a $\psi$ to which $u$
was close at time $0$, and then assumed that $u$ remained close to the
soliton evolution of $\psi$ over time.  This required us to make use of
the fact that $I\psi$ is close to $\psi$, which in turn forced us to require
the condition $\sigma \gtrsim N^{-C}$ for some large $C$.  We must
eliminate this condition in order to obtain the full force of the
theorem.  We will therefore find a $\psi^t$ which is close to $u$ for each
$t$, and study the equation by which this $\psi_t$ moves in time.  Define
$\psi_0(x)$ to be the standard ground state solution to equation (\ref{ODE})
centered at $0$.  

We begin by restating \ref{H1} in a form which will
be more convenient:
\begin{lemma}[Weinstein, \cite{W}] \label{lemmaW}
Let $\psi \in \Sigma$, and let $w \in H^1$ such that
$\|w\|_{H^1} \ll 1$ and $\langle w, (\psi^2)_x \rangle = 0$.  Then 
$$\Ly({\psi} + w) - \Ly(\psi_0) = \Ly({\psi} + w) - \Ly({\psi})
\sim \|w\|^2_{H^1}.$$
\end{lemma}
We will use the next lemma to find an appropriate ground state
${\psi}$ for each $t$ such that $u$ is close to ${\psi}$
and $w = u - {\psi}$ satisfies an appropriate orthogonality
condition.  Note that, since we will be studying $Iw$, not $w$, we will
require $\langle Iw, (\psi^2)_x \rangle = 0$ instead of $\langle w, (\psi^2)_x \rangle = 0$.
\begin{lemma}\label{lemmax0}
Let $u \in H^s$, and suppose $\mathrm{dist}_{H^s} (u, \Sigma) \ll N^{s
- 1}$ with $N$ sufficiently large.  Then $ u = {\psi} + w$ where
${\psi}$ is a ground state, $\langle w, I ({\psi}^2)_x
\rangle = 0$ and $\|Iw\|_{H^1} \lesssim N^{1 - s} \mathrm{dist}_{H^s}
(u , \Sigma) \ll 1$.
\end{lemma}

\begin{Proof}
Define $d(u,v) = \|I(u - v)\|_{H^1}$.  Then $d(u, \Sigma) \leq N^{1 -
s} \mathrm{dist}_{H^s} (u , \Sigma) \ll 1$.  So, as in \cite{CKSTT3},
there exists a $\psi'$ which minimizes $d(u, \Sigma)$.  By the
translation invariance of the problem, we may assume that this minimum
occurs at $\psi_0$.  Note that the tangent space to $\Sigma$ at
$\psi_0$ is spanned by $\psi_{0,x}$.  Therefore, if we differentiate
$$d(u, \psi)^2 = \langle I(u - \psi, I(u - \psi) \rangle +
\langle \p_x I(u - \psi), \p_x I (u - \psi) \rangle $$
in the $\psi_{0,x}$ direction, we will get $0$:
$$ \langle I (u - \psi_0) , I \psi_{0,x} \rangle + \langle \p_x I (u -
\psi_0), \p_x I \psi_{0,x} \rangle = 0 . $$
Let $\tilde{w} = u - \psi_0$.  Then, since $\psi_0 - \psi_{0,xx} - \psi_0^2 = 0
$, after integration by parts we get
$$\langle \tilde{w}, I^2 (\psi_0^2)_x \rangle = 0.$$
This is almost what we want; we would like to replace the $I^2$ in the
above equation by $I$.  To do so, we will perturb $\psi_0$ slightly.
Write ${\psi} = \psi_0(x - x_0)$, $w = u - {\psi}$ and $q =
{\psi} - \psi_0$.  
We want to solve $\langle w, I({\psi}^2)_x \rangle = 0$.  Using
what we know---$\langle \tilde{w}, I^2(\psi_0^2)_x \rangle = 0$---and
some algebra, what we want to solve for is:
$$ \langle q , I ((\psi_0 + q)^2)_x \rangle = \langle \tilde{w} , I ((\psi_0
+ q)^2)_x - I (\psi_0^2)_x \rangle + \langle I \tilde{w}, (\psi_0^2)_x -
I(\psi_0^2)_x \rangle . $$
Note that the last term is ${\cal O} (N^{-100} \|I\tilde{w} \|_{H^1} )=
{\cal O} (N^{-99} \mathrm{dist}_{H^s} (u , \Sigma))$, 
because $I - 1$ is almost the identity on $\psi_0$.  For the left-hand
side, note that $q = {\psi} - \psi_0 = -x_0 \psi_{0,x} + {\cal
O}_{H^2} (|x_0|^2)$, where ${\cal O}_{H^2}$ denotes the order of the
$H^2$ norm of a function.  Moreover, $((\psi_0 + q)^2)_x -(\psi_0^2)_x =
2x_0 (\psi_x^2 + \psi \psi_{xx}) + {\cal O}_{H^2} (|x_0|^2)$.
Therefore, the equation we wish to solve is
$$\langle x_0 \psi_{0,x} + {\cal O}_{H^2} (|x_0|^2) , I((\psi_0 + q)^2)_x
\rangle - x_0 \langle \tilde{w} , I(2 \psi_{0,x}^2 + \psi_0 \psi_{0,xx} + {\cal
O}_{H^2} (|x_0|^2)) \rangle = {\cal O} (N^{-99} \mathrm{dist}_{H^s} (u
, \Sigma)) . $$
Since $\psi_{0,x}^2 + \psi_0 \psi_{0,xx}$ is Schwartz, $$\langle I \tilde{w} , 2
\psi^2_{0,x} + \psi_0 \psi_{0,xx} \rangle = {\cal O} (\|I\tilde{w}\|_{H^1}) ={\cal O}
(d(u, \Sigma)) \ll 1.$$
On the other hand, $$\langle \psi_{0,x}, I (\psi_0^2)_x \rangle \sim \|\psi_0
\|^2_{W^{2,4}} + {\cal O}(N^{-100}),$$
which is an absolute constant that is not close to zero. 
So in the end, we get
$$x_0 (\langle \psi_{0,x}, I (\psi_0^2)_x \rangle - \langle I \tilde{w} , 2
\psi^2_{0,x} + \psi_0 \psi_{0,xx} \rangle ) = {\cal O} (N^{-99}
\mathrm{dist}_{H^s} (u , \Sigma) ) + {\cal O} (|x_0|^2), $$
where the coefficient of $x_0$ on the left-hand side is close to a
constant independent of $\tilde{w}$.
Therefore, by the inverse function theorem, we find that there is an
$x_0 \sim {\cal O} (N^{-99} \mathrm{dist}_{H^s} (u , \Sigma) )$ which
solves this equation, and then since $\|\psi - \psi _0\|_{H^2} = {\cal O}
(N^{-99} \mathrm{dist}_{H^s} (u , \Sigma) )$ the functions
${\psi} = \psi_0 ( x - x_0)$ and $w = u - {\psi}$ will
satisfy all the desired conditions.
\end{Proof}

We apply this lemma at each time $t$ such that
$\mathrm{dist}_{H^s} (u , \Sigma) \ll N^{s - 1}$ to write 
$u(x, t)~=~{\psi}^t(x)~+~w(x,t)$.   We will redefine $Q(x, t)$ by:
$$u(x,t) = Q(x,t) + w(x,t) = \psi_0(x - t - x_0(t)) + w(x,t)$$
For this section, we will redefine $E_N(t)$ in order to eliminate our
dependence on the closeness of $\psi$ and $I\psi$ and to reflect the more
precisely chosen error function $w(t)$ found in the above lemma.  We
therefore set\footnote{Compare to (\ref{E}).}
\begin{equation}
E_N(t) = \Ly(Q(t) + Iw(t)) . \label{E2}
\end{equation}
Note that, by (\ref{lemmax0}), for each $t$ such that
$\mathrm{dist}_{H^s} (u , \Sigma) 
\ll N^{s - 1}$,  $\langle Iw(t), (Q^2(t))_x \rangle = 0$ and
$\|Iw\|_{H^1} \ll 1$.  Therefore, by \ref{lemmaW}, $|E_N(t) - L(Q(t) )| \sim
\|Iw\|_{H^1}^2$.  In particular, at $t=0$, we have 
$$|E_N(0) - \Ly(Q(0))| \sim \|Iw\|^2_{H^1} \lesssim N^{2 - 2s}
\sigma^2$$
To prove the theorem, we will need the following lemma, a refinement
of Lemmas \ref{lemma1} and \ref{lemma2}:
\begin{lemma}\label{lemma3}
Suppose there is a $t_0 \in \R$ and a $\tilde{\sigma}$ with $0 < \tilde{\sigma} \ll 1$ such
that $|E_N(t_0) - \Ly(\psi_0)| \lesssim  \tilde{\sigma}^2$. Then there
exists a $\delta > 0$ depending only on $s$ such that 
$$E_N(t_0 + \delta) - E_N(t_0) = {\cal O} (\frac1{N^{1 - \epsilon}}
\tilde{\sigma}^2).$$
\end{lemma}
We will assume this lemma for now and conclude the proof of the
theorem:

\medskip
\begin{Proof} (of Theorem \ref{main2})
Once again, we set $\tilde{\sigma} = N^{1 - s} \sigma$.
As in the proofs of Propositions \ref{prop1} and \ref{prop2}, we can
iterate the result of Lemma \ref{lemma3}.  In this case, for $Q(x,t)=
\psi_0(x - t -x_0(t))$,
we obtain
that $|E_N(t) - \Ly(Q(t))| \lesssim N^{2 - 2s} \sigma^2$ for all $t$ such
that $t \ll N^{1 - \epsilon}$.  So, by Lemma~\ref{lemmaW}, 
$$\mathrm{dist}_{H^s} (u , \Sigma) \lesssim \|w\|_{H^s} \lesssim
N^{1-s} \sigma , $$ for all $t \ll N^{1 - \epsilon}$.  
We therefore can optimize for $N$ under only the two conditions:
\begin{align}\label{conditions3}
t \ll N^{1 - \epsilon} \ \ \ \ \ \ \ \ \ & \ \ \ \ \ \ \ \ \ N^{2 - 2s} \sigma ^2 \ll 1 .
\end{align}
Contrast these conditions with (\ref{conditions1}) and (\ref{conditions2}).
Note that we have now eliminated the condition $\sigma \ll N^{-C}$
and therefore we obtain $$\mathrm{dist}_{H^s} (u , \Sigma) \lesssim t^{1
- s + \epsilon} \sigma, $$ for all $t \ll \sigma^\frac1{1 - s -
\epsilon}$, as claimed. 
\end{Proof}

It thus remains only to prove Lemma~\ref{lemma3}:

\begin{Proof}(of Lemma~\ref{lemma3})
We write $Q(x, t) = \psi_0(x - t - x_0(t))$ and $w(x, t) = u(x, t) -
Q(x,t)$.  
Then $w(t)$ satisfies the difference equation:
\begin{equation}
w_t + w_{xxx} + \p_x(w (w + 2 Q) ) + \dot{x_0} Q_x = 0
\label{diffeqn2}\footnote{We use the notation $\dot{x_0}$ to mean the
ordinary derivative $\frac{dx}{dt}$.} 
\end{equation}
We know that $\|Iw(t_0)\|_{H^1} \lesssim \tilde{\sigma} = N^{1 -
s} \sigma$.  As before, we start by proving that the $X^{1,\frac12 +
\epsilon}$ norm of $Iw$ is controlled.  
\begin{claim} 
$$\|Iw\|_{X^{1, \frac12 +\epsilon}_{[t_0 - \delta, t_0 + \delta]}}
\lesssim \tilde{\sigma}. $$
\end{claim}
\begin{Proof}
As in each of the two previous claims, we use the standard $X^{s, b}$
estimates to obtain:
\begin{align*}
\|Iw\|_{X^{1, \frac12 +\epsilon}_{[t_0 - \delta, t_0 + \delta]}}
& \lesssim \|Iw(t_0)\|_{H^1} + \delta^{\epsilon} \|I(w_t +
w_{xxx})\|_{X^{1, -\frac12 +2 \epsilon}_{[t_0 - \delta, t_0 +
\delta]}}\\
& \lesssim \tilde{\sigma} + \delta^{\epsilon} \|I(w(w + 2
Q))_x\|_{X^{1, - \frac12 + 2 \epsilon}_{[t_0 - \delta, t_0 +
\delta]}} + \delta^{\epsilon} \|\dot{x}_0 (t) Q_x\|_{X^{1, -
\frac12 +2 \epsilon}_{[t_0 - \delta, t_0 + \delta]}} .
\end{align*}
Note that the first term on the right-hand side is the same as in
Claim~\ref{claim2} and can be estimated in exactly the same way.  For
the second term, we will prove that for each $t \in [t_0 - \delta, t_0
+ \delta]$, $\| \dot{x}_0 (t) \|
\lesssim \|Iw(t)\|_{H^1_x}$.  Then we will have
\begin{align*}
\|\dot{x}_0 (t) Q_x \|_{X^{1, - \frac12 + 2 \epsilon}_{[t_0 -
\delta, t_0 + \delta]}} & \lesssim \| \ \|Iw(t)\|_{H^1} Q_x \|_{X^{1,
- \frac12 + 2 \epsilon}_{[t_0 - \delta, t_0 + \delta]}} \\
& = \left \| \left ( \frac{ \langle \xi
\rangle}{\langle \tau - \xi ^3 \rangle^{\frac12 - 2 \epsilon}}
\left ( \widehat{\|Iw\|_{H^1_x}} \ast_\tau \tilde{Q}_x(\xi) \right ) (\tau) \right )
\right \|_{L^2_\tau L^2_\xi}\\
& \leq \| \ \widehat{\|Iw\|_{H^1_x}}\|_{L^1_\tau} \ \ \|\frac{ \langle \xi
\rangle}{\langle \tau - \xi ^3 \rangle^{\frac12- 2 \epsilon}}
\tilde{Q}_x(\xi, \tau - a) \|_{L^\infty_a L^2_{\xi, \tau}}\\
& \leq  C \|Iw\|_{L^\infty_{t,[t_0 - \delta, t_0 +
\delta]}  H^1_x} \\
& \leq C \|Iw\|_{X^{1, \frac12 + \epsilon}_{[t_0 - \delta, t_0 +
\delta]}} .
\end{align*}
The third line makes use of Minkowski's inequality for integrals, and
the fourth takes advantage of the fact that $Q(x,t)$ and all of its
$x$-translates are uniformly bounded in $X^{1, - \frac12 + 2 \epsilon}$ space.  The last
step is due to the standard estimate $\|f\|_{L^\infty_tH^1_x} \lesssim 
\|f\|_{X^{1, \frac12 + \epsilon}}$.  This argument allows us to
conclude that: 
$$\|Iw\|_{X^{1, \frac12 +\epsilon}_{[t_0 - \delta, t_0 + \delta]}}
\lesssim \tilde{\sigma} + \delta^{\epsilon} \|I(w(w + 2
Q))_x\|_{X^{1, - \frac12 + 2 \epsilon}_{[t_0 - \delta, t_0 +
\delta]}} + C \delta^{\epsilon} \|Iw\|_{X^{1, \frac12 +
\epsilon}_{[t_0 - \delta, t_0 + \delta]}}, $$
and we can then complete the proof of the claim via a continuity
argument.   
Therefore, to check that $\|Iw\|_{X^{1, \frac12 + \epsilon}_{[t_0 -
\delta, t_0 + \delta]}} \lesssim \tilde{\sigma}$ , we need only prove
that, for each~$t$, $|\dot{x}_0(t)| \lesssim \|Iw(t)\|_{H^1_x}$.

To do so, write $$\theta (x,t) = w(x + t + x_0(t),t) = u(x + t
+ x_0(t), t) - Q(x + t + x_0(t), t) = u(x + t + x_0(t), t) - \psi_0(x)$$
Then $\theta$ satisfies:
$$\theta_t + \theta_{xxx} + (\theta(\theta + 2 Q))_x =
\dot{x}_0(t)u_x + \theta_x . $$
Recall that $w$ satisfies $\langle w, I(Q^2)_x \rangle = 0$.
Differentiating in time, we see that, for each $t$, 
$$\langle \theta_t, I(\psi^2_0)_x \rangle = 0$$
Plugging in for $\theta_t$ and simplifying, we obtain: 
$$\dot{x}_0(t) \langle u_x, I(\psi_0^2)_x \rangle = \langle
\theta_{xxx}, I(\psi_0^2) _x \rangle + \langle (\theta(\theta + 2
\psi_0))_x , I (\psi_0^2)_x \rangle + \langle \theta_x, I(\psi_0^2)_x
\rangle , $$
i.e.
$$\dot{x}_0(t) = \frac1{\langle \psi_{0,x} + \theta_x, I (\psi_0^2)_x
\rangle} \left (\langle I\theta_x, (\psi_0^2)_{xxx} + (\theta +
2\psi_0)(\psi_0^2)_x + (\psi_0^2)_x \rangle + \langle \theta, I (\theta
+ 2 \psi_0)_x (\psi_0^2)_x \rangle \right ) . $$
Note that the numerator is controlled by $\|I \theta\|_{H^1}$ and
that the denominator is of a size greater than an absolute constant.
Therefore, we conclude 
that $\dot{x}_0(t)$ is indeed controlled by $\|I\theta(t)\|_{H^1} =
\|Iw(t)\|_{H^1}$ as claimed.  \end{Proof}

The final step in the proof of the lemma is to take a $\delta$ step
forward in $t$.  We want to prove
$$
E_N(t_0 + \delta) - E_N(t_0) = {\cal O} (\frac1{N^{1 - \epsilon}}
\tilde{\sigma}^2) .
$$
Recall that $E_N(t) = \Ly(Q + Iw) = \Ly (Q(x+t,t) + Iw(x+t,t))$.  
Also recall that $\Omega(f)(t) = \p_t (\Ly(f)(t)) = 2 \langle f_t, f -
f_{xx} - f^2 \rangle$.  
Therefore, 
\begin{align*}
\p_t E_N(t) & = \Omega((Q + Iw)(x + t,t)) \\
& = \langle \p_t ((Q + Iw)(x +t, t)), Q + Iw - Q_{xx}
-Iw_{xx} -2QIw -(Iw)^2 -Q^2\rangle \\
& = 2 \langle -\dot{x_0}(t) Q_x + I(-w_{xxx} - (w(w + 2Q))_x +
\dot{x_0}(t) Q_x), Iw - Iw_{xx} - (Iw)(2Q + Iw) \rangle \\
& = 2 \langle \dot{x_0}(t) ( IQ_x - Q_x), (Iw - Iw_{xx} -(Iw)(2Q + Iw)
\rangle + \\
& \phantom{moveit!} \langle I(w_x - w_{xxx} - (w(w +
2Q))_x, Iw - Iw_{xx} -(Iw)(2Q + Iw) \rangle .
\end{align*}
By integration by parts and the fact that $I$ is almost the identity
on $Q$, the first term (when integrated in $t$)
will be controlled by $C N^{-100} \|Iw\|^2_{X^{1,
\frac12 + \epsilon}_{[t_0 - \delta, t_0 + \delta]}}$ and will therefore be fine for our estimates.  

Note also that the second term is a polynomial of degree at least $2$
in $w$.  Therefore, as in Section~\ref{section2}, we will be done if we
can prove that for all $\gamma$ such that $\|\gamma\|_{X^{1,
\frac12 + \epsilon}_{[t_0 - \delta, t_0 + \delta]}} \leq 1$, 
$$\int_{t_0}^{t_0 + \delta} \langle I(\gamma_x - \gamma_{xxx} -
(\gamma(\gamma + 
2Q))_x), I\gamma- I\gamma_{xx} -(I\gamma)(2Q + I\gamma)
\rangle \lesssim \frac1{N^{1 - \epsilon}} . $$ 

To do so, let $v = \gamma + Q$.  Then 
$$\gamma_x - \gamma_{xxx}-(\gamma(\gamma + 2 q))_x = v_x - v_{xxx} - (v^2)_x$$ 
and
$$I\gamma - I\gamma_{xx} - (I\gamma)(2 Q + I\gamma) = Iv - Iv_{xx} -
(Iv)^2 + 2(Iv)(IQ - Q) - (IQ -Q)^2 -(IQ^2 - Q^2) . $$ 
Therefore 
$$\hskip-2.1in
\int_{t_0}^{t_0 + \delta} \langle I(\gamma_x - \gamma_{xxx} -
(\gamma(\gamma + 
2Q))_x, I(\gamma - \gamma_{xx} -\gamma(2Q + I\gamma{w}))
\rangle \ dt = $$
$$ \hskip0.3in = \int_{t_0}^{t_0 + \delta} \langle I(v_x - v_{xxx}-(v^2)_x), \left ( (Iv -
Iv_{xx} - (Iv)^2) + 2 (Iv)(I - 1)Q + ((I - 1)Q)^2  + (I - 1) Q^2 \right
) \rangle \ dt . 
$$
Then, once again, because $I - 1$ is nearly $0$ on $Q$, the second
third and fourth terms are controlled.  For the remaining term, note that, by
integration by parts $\langle Iv_x, Iv - Iv_{xx} - (Iv)^2) \rangle$ is
zero, so the last term to be estimated is 
\begin{equation}\label{last}
\int_{t_0}^{t_0 + \delta} \langle I(- v_{xxx}-(v^2)_x), Iv -
Iv_{xx} - (Iv)^2) \rangle dt . 
\end{equation}
But this is exactly the quantity estimated in Lemma~\ref{lemma1}.
Recall that the multilinear estimates used to prove those estimates did
not depend 
on the properties of the function $u$ except that $\| u\|_{X^{1,
\frac12 + \epsilon}_{[t_0 - \delta, t_0 + \delta]}} \lesssim 1$.  The
conclusion was that   
$$\int_{t_0}^{t_0 + \delta} \langle I(- v_{xxx}-(v^2)_x), Iv -
Iv_{xx} - (Iv)^2) \rangle dt = {\cal O} (\frac1{N^{1 - \epsilon}}).$$
Since $\| v\|_{X^{1,
\frac12 + \epsilon}_{[t_0 - \delta, t_0 + \delta]}}$ is indeed
controlled by a constant, by the estimates in the proof of
Lemma~\ref{lemma1} the 
quantity~(\ref{last}) is also controlled by $\frac1{N^{1 - \epsilon}}$.
This concludes the proof of Lemma~\ref{lemma3} and, at last, the main
theorem. \end{Proof}

\bibliographystyle{amsalpha}

\end{document}